\documentclass[12pt]{article}
\usepackage{amsmath,amssymb, amsthm}
\newtheorem{theorem}{Theorem}[section]

\newtheorem{lemma}[theorem]{Lemma}

\numberwithin{equation}{section}

\def\P{{{\mathcal{P}}}}
\def\S{{{\mathcal{S}}}}

\def\E_1{{{\mathcal{P}_1}}}
\def\E{{{\mathcal{E}}}}

\begin{document}
\title{Incidences between points and generalized spheres over finite fields and related problems}

\author{
Nguyen Duy Phuong\thanks{ This research was supported by Vietnam National Foundation for Science and Technology Development grant
    Email: {\tt duyphuong@vnu.edu.vn
}}
\and 
    Thang Pham\thanks{EPFL, Lausanne. Research partially supported by Swiss National Science Foundation
Grants 200020-144531 and 200021-137574.
    Email: {\tt thang.pham@epfl.ch}}
  \and
    Le Anh Vinh\thanks{ This research was supported by Vietnam National Foundation for Science and Technology Development grant.
    Email: {\tt vinhla@vnu.edu.vn, vinh@math.harvard.edu
}}}
\date{}
\maketitle

\noindent

\begin{abstract}
\noindent
Let $\mathbb{F}_q$ be a finite field of $q$ elements where $q$ is a large odd prime power and  $Q =a_1 x_1^{c_1}+\cdots+a_dx_d^{c_d}\in \mathbb{F}_q[x_1,\ldots,x_d]$, where $2\le c_i\le N$,   $\gcd(c_i,q)=1$, and  $a_i\in \mathbb{F}_q$ for all $1\le i\le d$. A \textit{$Q$-sphere} is a set of the form $\left\lbrace x\in \mathbb{F}_q^d~|~ Q(x-b)=r\right\rbrace$, where  $b\in \mathbb{F}_q^d, r\in \mathbb{F}_q$. We prove bounds on  the number of incidences between a point set $\P$ and a $Q$-sphere set $\S$, denoted by $I(\P,\S)$, as the following. 
  $$\left\vert I(\P,\S)-\frac{|\P||\S|}{q}\right\vert\le q^{d/2}\sqrt{|\P||\S|}.$$ 
  We prove this estimate by studying the spectra of directed graphs. 
We also give a version of this estimate over finite rings $\mathbb{Z}_q$ where $q$ is an odd integer. As a consequence of the above bounds, we give an estimate for the pinned distance problem. In Sections $4$ and $5$,  we  prove a bound on the number of incidences between a random point set and a random $Q$-sphere set in $\mathbb{F}_q^d$. We also study the finite field analogues of some combinatorial geometry problems, namely, the number of generalized isosceles triangles, and the existence of a large subset without repeated generalized distances.  
\end{abstract}
\section{Introduction}

Let $\mathbb{F}_q$ be a finite field of $q$ elements where $q$ is a large odd prime power. Let $P$ be a set of points, $L$  a set of lines over $\mathbb{F}_q^d$, and $I(P,L)$  the number of incidences between $P$ and $L$. Bourgain, Katz, and Tao \cite{bourgain-katz-tao} proved that for any $0<\alpha<2$ and $|P|,|L|\le N=q^{\alpha}$, $I(P,L)\lesssim N^{3/2-\epsilon}$, where $\epsilon=\epsilon(\alpha)$. By employing the Erd\H{o}s-R\'enyi graph (see $2.1$ for the definition),  the third author \cite{vinh-incidence} improved this bound in the case $1\le \alpha\le 2$, and gave the following estimate.
\begin{theorem}\label{m1}
Let $\mathcal{P}$ be a set of points and $\mathcal{L}$ a set of lines in
  $\mathbb{F}^2_q$. Then we have
  \begin{equation}
    I(\mathcal{P}, \mathcal{L}) \leq \frac{|\mathcal{P}||\mathcal{L}|}{q} +
    q^{1/2}\sqrt{|\mathcal{P}| |\mathcal{L}|} \nonumber
  \end{equation}
\end{theorem}
The above result was also proved for  points and hyperplanes, and for points and $k$-subspaces (see \cite{threepoints,vinh-incidence} for more details). 

Let $Q =a_1 x_1^{c_1}+\cdots+a_dx_d^{c_d}\in \mathbb{F}_q[x_1,\ldots,x_d]$, where $2\le c_i\le N$, for some constant $N>0$,  $\gcd(c_i,q)=1$, and  $a_i\in \mathbb{F}_q$ for all $1\le i\le d$. We define the \textit{generalized sphere}, or \textit{$Q$-sphere}, centered at $b=(b_1,\ldots,b_d)$ of radius $r\in \mathbb{F}_q$ to be the set $\lbrace x\in \mathbb{F}_q^d~|~Q(x-b)=r\rbrace .$
The main purpose of this paper is to give a similar bound on the number of incidences between points and generalized spheres by employing the spectral graph  method. With the same method,  we also consider some related problems in Sections $4$ and $5$. Our main result is the following.
\begin{theorem}\label{dinhly1}
Let $\P$ be a set of points and $\S$ a set of $Q$-spheres with arbitrary radii in  $\mathbb{F}_q^d$.  Then the number of incidences between points and spheres satisfies  \begin{equation}\label{mot}\left\vert I(\P,\S)- \frac{|\P||\S|}{q}\right\vert \le q^{d/2}\sqrt{|\P||\S|} .\end{equation}
\end{theorem}
In the case $Q(x)=\sum_{i=1}^d x_i^2$, Cilleruelo et al. \cite{iosevichetal} have independently proved (\ref{mot}). In this case, we also obtain a similar estimate over finite rings (see \cite{pham} for the Szemer\'{e}di-Trotter theorem over finite rings).
\begin{theorem}\label{dinhly2}
Let $\mathcal{P}$ be a set of points and $\mathcal{S}$ a set of spheres with arbitrary radii in $\mathbb{Z}_q^d$, $q$ is an odd integer. Then the number of incidences between points and spheres satisfies  \[\left\vert I(\mathcal{P}, \mathcal{S})- \frac{|\P||\S|}{q}\right\vert\le \sqrt{2\tau(q)}\frac{q^d}{\gamma(q)^{d/2}}\sqrt{|\P||\S|},\] where $\gamma(q)$ is the smallest prime divisor of $q$, and $\tau(q)$ the number of divisors of $q$.
\end{theorem}
\paragraph{Generalized pinned distances:} Let $P(x)\in\mathbb{F}_q[x_1,\ldots,x_d]$ be a polynomial and $\E\subset \mathbb{F}_q^d$. Given $x\in \mathbb{F}_q^d$, we denote the pinned $P$-distance set  determined by $\E$ and $x$ by $$\Delta_P(\E,x)=\lbrace P(y-x)\in\mathbb{F}_q ~|~y\in \E\rbrace.$$ We are interested in finding the elements $x\in \mathbb{F}_q^d$ and the size of $\E\subset \mathbb{F}_q^d$ such that $\Delta_P(\E,x)\gtrsim q$. In the case $P(x)=x_1^2+\cdots+x_d^2$, Chapman   et al. \cite{chap} proved that for any subset $\E\subset \mathbb{F}_q^d$ such that $|\E|\ge q^{(d+1)/2}$,  there exists a subset $\E'\subset \E$ such that $|\E'|\sim |\E|$, and  for every $y\in \E'$ we have $|\Delta_P(\E,y)|>\frac{q}{2}$. Cilleruelo et al. \cite{iosevichetal} reproved the same result  using their bound on number of incidences between points and spheres.

 In this general setting, the main difficulty in this problem is that we do not know the explicit form of the polynomial $P(x)$. Koh and Shen \cite{kod} found some conditions on $P(x)$  to obtain the desired bound. We remark that if $P$ is a diagonal polynomial of the form $\sum_{j=1}^d a_jx_j^s$, the conditions of Koh and Shen  are satisfied. However, if we consider the polynomial $Q(x)=\sum_{j=1}^d a_jx_j^{c_j}$, where the exponents $c_j$ are distinct, then we have not found any reference which shows that those conditions are satisfied.

As a consequence of Theorem \ref{dinhly1}, the following result can be derived in a similar way to how \cite{iosevichetal} derived their result from their bound on the number of  incidences between points and spheres. It generalizes the pinned distance results of \cite{chap}.

\begin{theorem}\label{dinhly3}
Let $\E\subset \mathbb{F}_q^d$ with $|\E|> \sqrt{(1-c^2)/c^4}\cdot q^{(d+1)/2}$ for some $0< c< 1$. Then  the number of points $p\in \E$ satisfying $|\Delta_Q(\E,p)|>(1-c)q$ is at least $(1-c)|\E|$.
\end{theorem}

\paragraph{Incidences between a random point set and a random $Q$-sphere set:} It follows from Theorem \ref{dinhly1} that if $\P$ is a set of points and $\S$ is a set of $Q$-spheres such that $|\P||\S|> q^{d+2}$, then there exists at least one incidence pair $(p,s)\in \P\times \S$ with $p\in s$. We  improve  the bound $q^{d+2}$ in the sense that for any $\alpha\in (0,1)$ it suffices to take $t\ge C_\alpha q$ randomly chosen points and spheres over $\mathbb{F}_q^d$ to guarantee that the probability of no incidences is exponentially small, namely $\alpha^t$, when $q$ is large enough. We remark that the ideas in this part are similar to the case between points and lines in \cite{vinhrandom}. More precisely, our result is the following.
\begin{theorem}\label{random1}
For any $\alpha> 0$, there exists an integer $q_0=q_0(\alpha)$ and a  number $C_\alpha >0$ with the following property. When a point set $\P$ and a $Q$-sphere set $\S$ where $|\P|=|\S|=t\ge C_\alpha q$ are chosen randomly in $\mathbb{F}_q^d$, the probability of $ \lbrace(p,s)\in \P\times \S:p\in s\rbrace=\emptyset$ is at most $\alpha^t$, provided that $q\ge q_0$. 
\end{theorem}
\paragraph{Generalized isosceles triangles:} Given a set $\E$  of $n$ points in $\mathbb{R}^2$, let $h(\E)$ be the number of isosceles triangles determined by $\E$. Define
 $h(n)=\min_{|\E|=n}h(\E)$. Pach and Tardos \cite{pachtardos} proved that $h(n)=O(n^{2.136})$. In the present paper, we consider the finite field version of this problem. Let us give some notation: A \textit{$Q$-isosceles triangle} at a vertex $x$ is a triple of distinct elements $(x,y,z)\in \mathbb{F}_q^d\times \mathbb{F}_q^d\times \mathbb{F}_q^d$ such that $Q(x-y)=Q(x-z)$. We will show that for any subset $\E$ in $\mathbb{F}_q^d$ such that its cardinality is large enough, the number of isosceles triangles determined by $\E$ is $(1+o(1))|\E|^3/q$.
\begin{theorem}\label{tamgiaccan1}
Given a set of $n$ points $\E$ in $\mathbb{F}_q^d$, $d\ge 2$.  If $|\E|\gg q^{\frac{2(d+1)}{3}}$, then the number of isosceles triangles determined by $\E$ is $(1+o(1))|\E|^3/q.$
\end{theorem}
Here and throughout, $X\gtrsim Y$ means that $X\ge CY$ for some constant $C$ and $X\gg Y$ means that $Y=o(X)$, where $X,Y$ are viewed as functions of the parameter $q$. 
\paragraph{Distinct distance subset:} Given a set $\E$ of $n$ points in $\mathbb{R}^2$, let $g(\E)$ be the maximal cardinality of a subset $U$ in $\E$ such that no distance determined by $U$ occurs twice. Define $g(n)=\min_{|\E|=n} g(\E).$  Charalambides \cite{chara} proved that $n^{1/3}/{(\log n)}\lesssim g(n)\lesssim n^{1/2}/(\log n)^{1/4}$, where the upper bound is obtained from the  Erd\H{o}s distinct distances problem (see \cite{colon, leftmann} for more details, earlier results, and results in higher dimensions). In this paper, we study the finite field analogue of this problem.

Given a set of $n$ points $\E\subset \mathbb{F}_q^d$,  a subset $U\subset \E$ is called a \textit{distinct  $Q$-distance subset} if there are no four distinct points $x,y,z,t\in U$ such that $Q(x-y)=Q(z-t)$. Using the same method that Thiele used in $\mathbb{R}^2$ (see \cite[p.191]{pach} for more details), we show that for any large enough set $\E$  in $\mathbb{F}_q^d$, there exists a distinct $Q$-distance subset of cardinality at least $Cq^{1/3}$, for some constant $C$. More precisely, we have the following estimate.
\begin{theorem}\label{khoangcachphanbiet}
Let $\E\subset \mathbb{F}_q^d $, $d\ge 2$, $|\E|\gg q^{2(d+1)/3}.$ If   $U_Q \subset \E$ is a maximal distinct $Q$-distance subset of $\E$, then  $ q^{1/3} \lesssim |U_Q|\lesssim q^{1/2}$. 
\end{theorem}
 \paragraph{About the work of Cilleruelo, Iosevich, Lund, Roche-Newton, and Rudnev:} After we finished a draft of this paper, we learned that Cilleruelo et al. \cite{iosevichetal} had independently obtained the same bound for the number of incidences between points and spheres in the case $Q(x-y)=\sum_{i=1}^d(x_i-y_i)^2$, using the elementary method introduced in \cite{cil}. 
\section{Spectra of graphs and digraphs}
\subsection{Pseudo-random graphs}
Let us recall some notions about $(n,d,\lambda)$-graphs from Alon and Spencer in \cite{as}. Given  an undirected graph $G$, let $\lambda_1(G) \geq \lambda_2(G) \geq \ldots \geq \lambda_n(G)$ be
the eigenvalues of its adjacency matrix. The quantity $\lambda (G) = \max
\{\lambda_2(G), - \lambda_n(G) \}$ is called the second eigenvalue of $G$. A graph $G
= (V, E)$ is called an $(n, d, \lambda)$-graph if it is $d$-regular, has $n$
vertices, and the second eigenvalue of $G$ is at most $\lambda$. It is well known (see \cite[Chapter 9]{as} for more details) that if $\lambda$ is much smaller than the degree $d$, then $G$ has certain random-like properties.  For two (not necessarily disjoint) subsets of vertices $U, W \subset V$, let $e (U, W)$ be the number of ordered pairs $(u, w)$ such that $u \in U$, $w \in W$, and $(u, w)$ is an edge of $G$. For a vertex $v$ of $G$, let $N (v)$ denote the set of vertices of $G$ adjacent to $v$ and let $d (v)$ denote its degree. Similarly, for a subset $U$ of the vertex set, let $N_U (v) = N(v) \cap U$ and $d_U (v) = |N_U (v) |$. We first recall the following  well-known lemma (see, for example, \cite[Corollary 9.2.5]{as}).

\begin{lemma}\label{tool2}
  Let $G = (V, E)$ be an $(n, d, \lambda)$-graph. For any two sets $B, C
  \subset V$, we have
  \[ \left| e (B, C) - \frac{d|B | |C|}{n} \right| \leq \lambda \sqrt{|B| |C|}. \]
\end{lemma}

Let $PG(q,d)$ denote the projective space of dimension $d-1$ over the finite field $\mathbb{F}_q$. Let $\mathcal{ER}(\mathbb{F}_q^d)$ denote the graph with vertex set $PG(q,d)$, and two vertices $\textbf{x},\textbf{y}$ are connected by an edge if $\textbf{x}\cdot \textbf{y}=0$. In the case $d=2$, this graph is called \textit{Erd\H{o}s-R\'enyi graph}. The third author used the spectrum of $\mathcal{ER}(\mathbb{F}_q^d)$ and Lemma \ref{tool2} to prove Theorem \ref{m1} (see \cite{vinh-incidence} for more details).

  In order to prove Theorem \ref{dinhly2}, we use the \textit{sum-product graph} defined as the following.  The vertex set of the sum-product graph $\mathcal{SP}(\mathbb{Z}_q^{d+1})$ is the set $V(\mathcal{SP}(\mathbb{Z}_q^{d+1})) = \mathbb{Z}_q\times \mathbb{Z}_q^d$. Two vertices $U=(a,\textbf{b})$ and $V=(c,\textbf{d}) \in V(\mathcal{SP}(\mathbb{Z}_q^{d+1}))$ are connected by an edge, $(U,V) \in E(\mathcal{SP}(\mathbb{Z}_q^{d+1}))$, if and only if $a+c=\textbf{b}\cdot\textbf{d}$. Our construction is similar to that of Solymosi in $\cite{solymosi}$. We  have the following lemma about the spectrum of the sum-product graph $\mathcal{SP}(\mathbb{Z}_q^{d+1})$ (see \cite[Lemma 4.1]{vinh-sumproductgraph} for the proof).
\begin{lemma}\label{2.6}
For any $d\geq 1$, the sum-product graph $\mathcal{SP}(\mathbb{Z}_q^{d+1})$ is a
\[\left(q^{d+1}, q^d, \sqrt{2\tau(q)}\frac{q^d}{\gamma(q)^{d/2}} \right)-\mbox{graph}.\]
\end{lemma}

However, it seems difficult to use the spectrum of an undirected graph to analyze the number of incidences between points and $Q$-spheres, where $Q(x)\in \mathbb{F}_q[x_1,\ldots,x_d]$ is an arbitrary diagonal polynomial. In the next subsection, we will introduce the Cayley graph and some notions from Vu \cite{van} to deal with this problem.

\subsection{Pseudo-random digraphs}
Let $G$ be a directed graph (digraph) on $n$ vertices where the in-degree and out-degree of each vertex are both $d$. The adjacency matrix  $A_G$ is defined as follows: $a_{ij}=1$ if there is a directed edge from $i$ to $j$, and zero otherwise. Let $\lambda_1(G),\ldots,\lambda_n(G)$ be the eigenvalues of $A_G$. These numbers are complex numbers, so we can not order them, but we have $|\lambda_i|\le d$ for any $1\le i\le n$. Define $\lambda_1(G)=d, \lambda(G):=\max_{|\lambda_i(G)|\ne d} |\lambda_i(G)|$. 

 A digraph $G$ is called a $(n,d,\lambda)$-digraph if it has $n$ vertices, the in-degree and out-degree of each vertex is $d$, and $\lambda(G)\le \lambda$.  

Let $G$ be a $(n,d, \lambda)$-digraph. For any two (not necessarily disjoint) subsets $U,W\subset V$, let $e(U,W)$ be the number of ordered pairs $(u,w)\in U\times W$ such that $\overrightarrow{uw}$ is an edge of $G$.  Vu \cite[Lemma 3.1]{van} developed a directed version of the Lemma \ref{tool2}.
\begin{lemma}\label{vanvu}
Let $G=(V,E)$ be a $(n,d,\lambda)$-digraph. For any two sets $B,C\subset V$, we have $$ \left \vert e(B,C)-\frac{d|B||C|}{n}\right\vert  \le \lambda \sqrt{|B||C|}. $$
\end{lemma}
Let $H$ be a finite abelian group and $S$ a subset of $H$. The {\it Cayley graph} is the digraph $C_S(H)=(H, E)$, where the vertex set is $H$, and  there is a directed edge from vertex $x$ to  vertex $y$ if and only if $y-x\in S$.   It is clear that every vertex of $C_S(H)$ has out-degree $|S|$.  We define the graph $C_Q(\mathbb{F}_q^{d+1})$ to be the Cayley graph with $H=\mathbb{F}_q\times \mathbb{F}_q^d$ and $S=\lbrace (x_0,x)\in \mathbb{F}_q\times \mathbb{F}_q^d~|~x_0+Q(x)=0\rbrace$, i.e. $$E(C_Q(\mathbb{F}_q^{d+1}))=\{((x_0,x), (y_0,y))\in H\times H ~| ~x_0-y_0+Q(x-y)=0\}.$$
We have the following result on the spectrum of  $C_Q(\mathbb{F}_q^{d+1})$. We reproduce the proof because this lemma is crucial to our main results.
\begin{lemma}(See \cite[Lemma 3.2]{vinherdos}.)\label{bodechinh}
For any odd prime power $q$, $d\ge 1$, then $C_Q(\mathbb{F}_q^{d+1})$ is a 
  \[(q^{d+1}, q^d, q^{d/2})-\mbox{digraph}.\]
\end{lemma}
With the same arguments, we obtain the following lemma for the graph we use in the proof of Theorem \ref{tamgiaccan1}.
\begin{lemma}\label{bodechinh11}
For any odd prime power $q$, $d\ge 1$, let $Q'(x_1,\ldots,x_{2d})$ be a polynomial in $\mathbb{F}_q[x_1,\ldots,x_{2d}]$ defined by $ Q'=Q(x_1,\ldots,x_d)-Q(x_{d+1},\ldots,x_{2d})$. Then $C_{Q'}(\mathbb{F}_q^{2d+1})$ is a \[(q^{2d+1}, q^{2d}, q^{d})-\mbox{digraph}.\]
\end{lemma}

\section{Proofs of Theorems \ref{dinhly1} and \ref{dinhly2}}\label{section2}
\paragraph{Proof of Theorem \ref{dinhly1}}
 We use the Cayley graph $C_Q(\mathbb{F}_q^{d+1})$ to prove Theorem \ref{dinhly1}. Let $\P=\lbrace(x_{i1},\ldots,x_{id})\rbrace_i$ be a set of $n$ points in $\mathbb{F}_q^d$, and $S=\lbrace(r_i, (y_{i1},\ldots,y_{id}))\rbrace_i$  a set of pairs of radii and centers representing $Q$-spheres in $\S$. Let $U=\left\lbrace(0,x_{i1},\ldots,x_{id})\right\rbrace_i \subset\mathbb{F}_q^{d+1}$ and $W=\left\lbrace(r_i, y_{i1},\ldots,y_{id})\right\rbrace_i \subset\mathbb{F}_q^{d+1}$. Then the number of incidences between points and $Q$-spheres is the number of edges between $U$ and  $W$ in $C_Q(\mathbb{F}_q^{d+1})$. Using Lemma \ref{vanvu} and \ref{bodechinh},  Theorem \ref{dinhly1} follows. 
\paragraph{Proof of Theorem \ref{dinhly2}}  We use the sum-product graph $\mathcal{SP}(\mathbb{Z}_q^{d+1})$ to prove Theorem \ref{dinhly2}.
We identify each point $(b_1,\ldots,b_d)$ in $\P$ with a vertex  $(-b_1^2-\cdots-b_d^2,b_1,\ldots,b_d)\in \mathbb{Z}_q^{d+1}$ of $\mathcal{SP}(\mathbb{Z}_q^{d+1})$, and each sphere $(x_1-a_1)^2+\cdots+(x_d-a_d)^2=r$ in $\S$ with  a vertex $(r-a_1^2-\cdots-a_d^2,-2a_1,\ldots,-2a_d)\in \mathbb{Z}_q^{d+1}$ of $\mathcal{SP}(\mathbb{Z}_q^{d+1})$. Let $U\subset\mathbb{Z}_q^{d+1}$ be the set of points corresponding to $\P$, and $W\subset\mathbb{Z}_q^{d+1}$ the set of points corresponding to $\S$. Then the number of incidences between points and spheres is the number of edges between $U$ and $W$ in the sum-product graph $\mathcal{SP}(\mathbb{Z}_q^{d+1})$. By Lemma \ref{tool2} and Lemma \ref{2.6}, Theorem \ref{dinhly2} follows.

\paragraph{Remark:} The authors have not found any reference for a version of Weil's theorem over finite rings $\mathbb{Z}_m^d$.  Therefore, it seems hard to prove Theorem \ref{dinhly1} for a more general polynomial $Q(x)$ over finite rings using directed graphs. We note that Lemmas \ref{bodechinh} and \ref{bodechinh11} also hold for the general case $Q(x_1,\ldots,x_d)=\sum_{i=1}^df_i(x_i)$, where $\deg(f_i)\ge 2, ~\gcd(\deg(f_i),q)=1$ for all $i$. Therefore, all of the results in this paper over finite fields also hold for these more.
\section{Generalized pinned distance problem}\label{section3}

\paragraph{Proof of Theorem \ref{dinhly3}:} First we prove that $$\frac{1}{|\E|}\sum_{p\in \E}|\Delta_Q(\E,p)|> (1-c^2)q.$$
We identify each point $p=(b_1,\ldots,b_d)\in \E $ with a point $(0,b_1,\ldots,b_d)\in \mathbb{F}_q^{d+1}$, and  each pair $(p=(b_1,\ldots,b_d),t)$ where $t \in \Delta_Q(\E,p)$ with  a point $(t, b_1,\ldots,b_d)\in \mathbb{F}_q^{d+1}$. Let $U\subset \mathbb{F}_q^{d+1}$ be the set of points  corresponding to $\E$, and $W\subset \mathbb{F}_q^{d+1}$ the set of points  corresponding to point-distance pairs. Then  $|U|=|\E|$, $|W|=\sum_{p\in \E}|\Delta_Q(\E,p)|$. Moreover, one can easily see that $U,W$ are vertex subsets of the Cayley digraph $C_Q(\mathbb{F}_q^{d+1})$. The number of edges between $U$ and $W$ is $|\E|^2$, since each point in $\E$ contributes $|\E|$ edges between $U$ and $W$.  It follows from Lemmas \ref{vanvu} and \ref{bodechinh} that 
\begin{eqnarray}\label{18}
|\E|^2 \le e(U,W)&\le& \frac{|U||W|}{q}+q^{d/2}\sqrt{|U||W|}.\nonumber\\
&=&\frac{|\E|\sum_{p\in \E}|\Delta_Q(\E,p)|}{q}+q^{d/2}\sqrt{|\E|\sum_{p\in \E}|\Delta_Q(\E,p)|}.
\end{eqnarray}
If $\frac{1}{|\E|}\sum_{p\in \E}|\Delta_Q(\E,p)|\le (1-c^2)q$, it follows from (\ref{18}) that 
\begin{eqnarray}
|\E|^2&\le& |\E|^2(1-c^2)+q^{(d+1)/2}|\E|\sqrt{(1-c^2)}\nonumber\\
|\E|&\le &\sqrt{\frac{(1-c^2)}{c^4}}q^{(d+1)/2}.\nonumber
\end{eqnarray}
This would be a contradiction. Therefore, \begin{equation}\label{hehehe}\sum_{p\in \E}|\Delta_Q(\E,p)|> (1-c^2)q|\E|.\end{equation} Let us define $\E':=\lbrace p\in \E: |\Delta_Q(\E,p)|>(1-c)q\rbrace$. Suppose  that $|\E'|<(1-c)|\E|$, so
\begin{equation}\label{pt11.3}
\sum_{p\in \E\setminus \E'}|\Delta_Q(\E,p)|\le (|\E|-|\E'|)(1-c)q,
\end{equation}
and \begin{equation}\label{pt11.2} \sum_{p\in \E'}|\Delta_Q(\E,p)|\le q|\E'|.\end{equation}
Putting (\ref{pt11.3}) and (\ref{pt11.2}) together, we obtain
$$\sum_{p\in \E}|\Delta_Q(\E,p)|\le  (1-c)q|\E|+cq|\E'|<(1-c)q|\E|+cq(1-c)|\E|=(1-c^2)q|\E|.$$
The theorem follows because this contradicts (\ref{hehehe}).
\section{Related Problems}

\subsection{Incidences between random points and $Q$-spheres}\label{section5}

To prove Theorem \ref{random1}, we need the following lemma (see \cite[Lemma 8]{hoi}, and \cite[Lemma 2.3]{vinhrandom} for more details).
\begin{lemma}\label{random}
Let $\lbrace G_n=G(U_n,V_n)\rbrace_{n=1}^\infty$ be a sequence of bipartite graphs with $|V_n|=|U_n|\to \infty$ as $n\to \infty$, and let $\bar{d}(G_n)$ be the average degree of $G_n$. Assume that for any $\epsilon>0$, there exists an integer $v(\epsilon)$ and a number $c(\epsilon)>0$ such that $$e(A,B)\ge c(\epsilon)|A||B|\frac{\bar{d}(G_n)}{|V_n|},$$  for all $|V_n|=|U_n|\ge v(\epsilon)$ and all $A\subset V_n, B\subset U_n$ satisfying $|A||B|\ge \epsilon |V_n|^2$. Then for any $\alpha >0$, there exist an integer $v(\alpha)$ and a number $C(\alpha)$ with the following property: if one chooses a random subset $S$ of $V_n$ of cardinality $t$ and a random subset $T$ of $U_n$ of the same cardinality $t$, then the probability of $G(S,T)$ being empty is at most $\alpha^t$ provided that $t\ge C(\alpha)|V_n|/\bar{d}(G_n)$ and $|V_n|\ge v(\alpha)$. 
\end{lemma}
We notice that the Lemma \ref{random}  also holds when  $\{G_n\}_n$ is a sequence of digraphs.
\paragraph{Proof of Theorem \ref{random1}:} Let $B_{q,d}$ be a bipartite digraph with vertex set $V(C_Q(\mathbb{F}_q^{d+1}))\times V(C_Q(\mathbb{F}_q^{d+1}))$, where $C_Q(\mathbb{F}_q^{d+1})$ is the Cayley graph defined as in Lemma \ref{bodechinh} and the edge set
\[\lbrace\left((x_0,x),(y_0,y)\right)\in \mathbb{F}_q^{d+1}\times \mathbb{F}_q^{d+1} ~|~(x_0-y_0)+Q(x-y)=0 \rbrace  .\] With the same identification of the point set and the $Q$-sphere set as in proof of Theorem \ref{dinhly1}, we obtain two corresponding sets $U$ and $W$, where $|U|=|\P|$, $|W|=|\S|$. Thus, the number of incidences between points and spheres is the number of edges between $U$ and $W$. By Lemma \ref{vanvu} and \ref{bodechinh}, we obtain
\begin{equation}\label{hello9}
 \left\vert e(U, W)- \frac{|U||W|}{q}\right\vert\le q^{d/2}\sqrt{|U||W|}.
\end{equation}
For any $\epsilon>0$ such that $|U||W|\ge \epsilon q^{2d+2}$ and $q^d\ge \frac{4}{\epsilon}$, we have from (\ref{hello9}) that 
$$e(U,W) \ge \frac{q^d}{2q^{d+1}}|U||W|=\frac{\bar{d}(B_{q,d})}{|V(B_{q,d})|}|U||W|.$$ Let $c(\epsilon)=1, v(\epsilon)\ge (\frac{4}{\epsilon})^{(d+1)/d}$, then the theorem follows from Lemma \ref{random}.

\subsection{Generalized isosceles triangles}

\paragraph{Proof of Theorem \ref{tamgiaccan1}:}

Let \[U=\{(1,x,x)\in 1\times \E\times \E\},~ W=\{(1,y,z)\in 1\times \E\times \E\}.\]  One can easily see that $|U|=|\E|,|W|=|\E|^2$. Let \[T_1=\{(1,x,x,1,y,z)\in 1\times \E \times \E \times 1\times \E \times \E : Q(x-y)=Q(x-z)\}.\] Then the cardinality of $T_1$ is the number of edges between the sets  $U$ and $W$ in the graph $C_{Q'}(\mathbb{F}_q^{2d+1})$ (defined as in Lemma \ref{bodechinh11}). It follows from Lemma \ref{vanvu} and \ref{bodechinh11} that \[\left\vert|T_1|-\frac{|U||W|}{q}\right\vert\le q^d\sqrt{|U||W|}.\] Thus, if $|\E|\gg q^{2(d+1)/3}$ then $|T_1|=(1+o(1))|\E|^3/q$. We notice that $T_1$ also contains the tuples $(1,x,x,1,x, y)$ with $Q(x-y)=0$ which correspond to the edges between the vertices $(1,x,x)\in U$ and $(1,x,y)\in W$. Let us denote the set of such tuples by $T_{err}$, then one can easily see that  $\frac{1}{2}|T_{err}|$  is the number of pairs $(x,y)\in \E\times \E$ such that $Q(x-y)=0$, since each pair $(x,y)$ with $Q(x-y)=0$ contributes two edges $((1,x,x),(1,x,y))$ and $((1,x,x),(1,y,x))$. It follows from Lemma \ref{vanvu} and \ref{bodechinh}  that \[\left\vert |T_{err}|-\frac{|\E|^2}{q}\right\vert\le q^{d/2}\sqrt{|\E|^2}.\] Thus, if $|\E|\gg  q^{2(d+1)/3}$ with $d\ge 2$, then $|T_{err}|=|\E|^2/q=o(1)|\E|^3/q$. Therefore, the number of $Q$-isosceles triangles determined by $\E$ is $(1+o(1))|\E|^3/q$.
\subsection{Distinct distance subset}
In order to prove Theorem \ref{khoangcachphanbiet}, we need the following theorem on the cardinality of a maximal independent set of a hypergraph due to Spencer \cite{spencer}.
\begin{theorem}\label{spencer}
Let $H$ be a $k$-uniform hypergraph with $n$ vertices and $m \ge n/k$ edges, and let $\alpha(H)$ denote the independence number of $H$. Then 
$$ \alpha(H) \ge \left(1-\frac{1}{k}\right) \left \lfloor \left(\frac{1}{k}\frac{n^k}{m}\right)^{\frac{1}{k-1}}\right \rfloor.$$
\end{theorem}
\paragraph{Proof of Theorem \ref{khoangcachphanbiet}:} Let \[T_2=\{(1, p_1, q_1, 1, p_2, q_2)\in 1\times \E \times \E \times 1 \times \E\times \E: Q(p_1-q_1)=Q(p_2-q_2)\}.\]  With the same arguments  in the proof of Theorem \ref{tamgiaccan1}, we obtain  $|T_2|\le \frac{|\E|^4}{q}+q^d|\E|^2$. Thus, if $|\E|\gg q^{(d+1)/2}$, then $$|T_2|=(1+o(1)) \frac{|\E|^4}{q}.$$ 
A $4$-tuple of distinct elements in $\E^4$ is called \textit{regular} if all six generalized distances determined are distinct. Otherwise, it is called \textit{singular}. Let $H$ be the $4$-uniform hypergraph on the vertex set $V(H)=\E$, whose edges are the singular $4$-tuples of $\E$. 

It follows from Theorem \ref{tamgiaccan1} that the number of $4$-tuples containing a triple induced an isosceles triangle  is at most $((1+o(1))|\E|^3/q)\cdot |\E|=(1+o(1))|\E|^4/q$ when $|\E|\gg q^{2(d+1)/3}$. Thus the number of edges of $H$ containing a triple induced an isosceles triangle is at most $(1+o(1))|\E|^4/q$.  On the other hand, since $T_2=(1+o(1))|\E|^4/q$ when $|\E|\gg q^{(d+1)/2}$, the number of $4$-tuples $(p_1,q_1,p_2,q_2)$ in $\E^4$ satisfying $Q(p_1-q_1)=Q(p_2-q_2)$ equals $(1+o(1))|\E|^4/q$ when $|\E|\gg q^{(d+1)/2}$. Thus, if $|\E|\gg q^{2(d+1)/3}$ with $d\ge 2$, then $$|E(H)|\le \frac{2|\E|^4}{q}.$$
It follows from Theorem \ref{spencer} that $$\alpha(H)\ge C \left \lgroup\frac{|\E|^4}{|E(H)|}\right\rgroup^{1/3}= Cq^{1/3},$$
for some positive constant $C$.
Since there is no repeated generalized distance determined by the independent set of $H$,  we have $|U_Q|\ge \alpha(H)\ge Cq^{1/3}$.

Moreover, it is easy to see that there is at least one repeated generalized distance determined by any set of  $\sqrt{2}q^{1/2}+1$ elements  since there are only $q=|\mathbb{F}_q|$ distances over $\mathbb{F}_q^d$. Thus, the theorem follows.

\bigskip
\section*{Acknowledgements.}
\thispagestyle{empty}
The authors would like to thank J\'{a}nos Pach and Frank de Zeeuw for many useful discussions and helpful comments. The authors are also grateful to the referee for useful comments and suggestions.

\end{document}